\documentclass[11pt]{amsart}

\usepackage{epigamath}

\usepackage[english]{babel}

\numberwithin{equation}{section}

\usepackage{enumitem}
\usepackage{mathdots}

\newtheorem{thm}{Theorem}[section]

\newtheorem{prop}[thm]{Proposition}

\theoremstyle{definition}

\newtheorem{defi}[thm]{Definition}

\theoremstyle{remark}

\newtheorem{rem}[thm]{Remark}

\newcommand{\C}{\mathbb{C}}

\newcommand{\M}{\mathcal{M}}
\newcommand{\E}{\mathcal{E}}

\newcommand{\CP}{\mathbb{C}\mathbf{P}}
\newcommand{\RP}{\mathbb{R}\mathbf{P}}

\renewcommand{\epsilon}{\varepsilon}
\newcommand{\R}{\mathbb{R}}
\newcommand{\Q}{\mathbb{Q}}
\newcommand{\Z}{\mathbb{Z}}
\newcommand{\ZZ}{\Z/2\Z}

\newcommand{\sC}{\mathscr{C}}
\newcommand{\sE}{\mathscr{E}}

\newcommand{\cA}{\mathcal{A}}
\newcommand{\cG}{\mathcal{G}}
\renewcommand{\geq}{\geqslant}
\renewcommand{\ge}{\geqslant}
\renewcommand{\leq}{\leqslant}
\renewcommand{\le}{\leqslant}

\renewcommand{\ss}{\mathrm{ss}}

\EpigaVolumeYear{6}{2022} \EpigaArticleNr{24} \ReceivedOn{December 2, 2021}
\InFinalFormOn{July 27, 2022}
\AcceptedOn{August 26, 2022}

\title{Maximality of moduli spaces of vector bundles on curves}
\titlemark{Maximality of moduli spaces of vector bundles}

\author{Erwan Brugall\'e}
\address{Universit\'e de Nantes, Laboratoire de
  Math\'ematiques Jean Leray, 2 rue de la Houssini\`ere, F-44322 Nantes Cedex 3, France}
\email{erwan.brugalle@math.cnrs.fr}

\author{Florent Schaffhauser}
\address{Mathematisches Institut, Universität Heidelberg, Im Neuenheimer Feld 205, 69120 Heidelberg, Germany}
\email{fschaffhauser@mathi.uni-heidelberg.de}

\authormark{E.~Brugall\'e and F.~Schaffhauser}

\AbstractInEnglish{We prove that moduli spaces of semistable vector bundles of coprime
rank and degree over a non-singular real projective curve are maximal
real algebraic varieties if and only if the base curve itself is
maximal. This provides a new family of maximal varieties, with members
of arbitrarily large dimension. We prove the result by comparing the
Betti numbers of the real locus to the Hodge numbers of the complex
locus and showing that moduli spaces of vector bundles over a maximal
curve actually satisfy a property which is stronger than maximality
and that we call Hodge-expressivity. 
We also give a brief account on 
other varieties for which this property was already known.}

\MSCclass{14P25, 14H60}

\KeyWords{Maximal real algebraic varieties, Hodge numbers, moduli spaces of vector bundles}

\acknowledgement{E.~Brugallé  is
partially supported by the grant TROPICOUNT of Région Pays de la
Loire, and the ANR project ENUMGEOM NR-18-CE40-0009-02. 

During the completion of this work, F. Schaffhauser has been partially supported by the University of Strasbourg Institute of Advanced Study, the Max Planck Institute for Mathematics in Bonn and DFG-funded SPP Project 1786 ``Homotopy theory and algebraic geometry''.}

\begin{document}



\maketitle

\begin{prelims}

\DisplayAbstractInEnglish

\bigskip

\DisplayKeyWords

\medskip

\DisplayMSCclass







\end{prelims}


\newpage

\setcounter{tocdepth}{1}

\tableofcontents


\section*{Notation}
A real algebraic variety is a pair $(X,\tau)$ where $X$ is a complex
algebraic variety equipped with an antiholomorphic involution
$\tau:X\to X$. We denote  by
$\R X$ the set of  its real points, that is to say the fixed point set
of $\tau$.
Throughout this note, complex algebraic varieties are endowed with
their Euclidean topology. Except otherwise stated, we consider cohomology with
coefficients in $\ZZ$. 
Given a smooth compact manifold $M$ and a
non-singular
projective complex  algebraic variety $X$, we denote by 
$P_{t}(M)$ the mod 2 Poincaré polynomial of $M$ and by
$H_{(x,y)}(X)$  the Hodge polynomial of $X$, \textit{i.e.}
\[
P_t(M)\ =\ \sum_{i\geq 0}b_i(M)t^i \quad\mbox{and}\quad
H_{(x,y)}(X)\ =\ \sum_{i,j \geq 0}h^{i,j}(X)x^iy^j.
\]

\section{Introduction}

\subsection{Maximal and Hodge-expressive real algebraic varieties}
Let $X$ be a real algebraic variety (not necessarily projective nor smooth). It is a consequence of Smith theory (see for example \cite[Section
  3.3]{Mang17}) that the mod 2 Betti numbers of $\R X$ and $X$ satisfy
the following
Smith--Thom inequality:
\[
\sum_{i= 0}^{\dim X}b_i(\R X)\ \leq\ \sum_{i= 0}^{2\,\dim X}b_i(X).
\]
When equality holds, the real algebraic variety $X$ is said to be
\textit{maximal}.
Maximal varieties constitute extremal objects in real algebraic
geometry that enjoy special properties. We refer to the book
\cite{Mang17} or to the survey
\cite{DK} for an account on the subject.

To date, not many maximal real algebraic varieties are
known. It is standard that Grassmannians equipped with their standard real structure
(\textit{e.g.} real projective spaces) are maximal, as well as 
 non-singular (or mildly singular) real toric varieties \cite{BFMV06,Fr22}.
One also easily constructs maximal real algebraic (hyperelliptic) curves of arbitrary
genus. Furthermore, the Jacobian of a real algebraic curve with
non-empty real part  is
maximal if and only if the curve is maximal \cite{GrHa83}.
Knowledge starts to be more fragmented in the case of surfaces,  we refer the
interested reader to \cite[Chapter 4]{Mang17} or
\cite[Section 3]{DK}, as well as to Section \ref{sec:discussion}. Things are of course
getting worse when increasing the dimension; the knowledge of the
authors essentially reduces to the classifications of maximal cubic
real hypersurfaces of dimension~3 and~4 \cite{Kra09,FinKha10}, and to the
 unpublished construction by
Itenberg and Viro, dating back to the 1990's, of maximal real projective hypersurfaces of arbitrary degree and
dimension.

\medskip
 The goal of this note is to  provide a new family of maximal
 real algebraic varieties with members of arbitrarily large
 dimension: moduli spaces of vector bundles of coprime rank and degree over a maximal
 real algebraic curve.

\begin{figure}[!h]
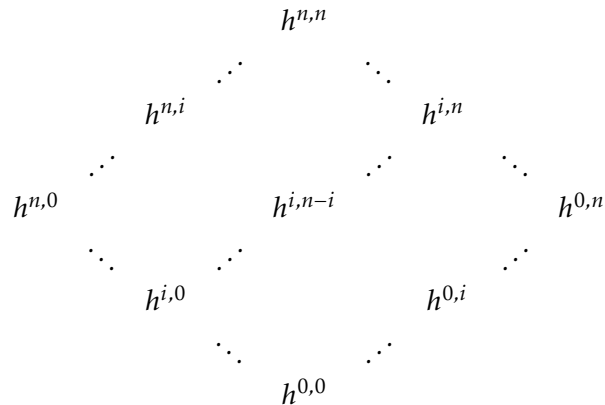

 \[
\begin{array}{ccccccccc}
&&&&h^{n,n}&&&& \\
&&&\iddots&& \ddots &&& \\
 &&h^{n,i} &&&&h^{i,n}\  &&  \\
 & \iddots & &&&\iddots&&\ddots&  \\
h^{n,0} &  &&& h^{i,n-i} & && & h^{0,n} \\
 &\ddots  &&\iddots &&  && \iddots &  \\
  & & h^{i,0} & && & h^{0,i} &  & \\
  & & & \ddots & & \iddots && & \\
   & & &  &  h^{0,0} &&& &
\end{array}
\]
\caption{Ascending diagonals of the Hodge diamond.} \label{Hodge_diamond}
\end{figure}

 Poincaré polynomials of the complex and real part of these moduli spaces have been computed in
\cite{AtBo83} and \cite{LiSc13} respectively (see also \cite{Baird_CJM}). 
  It turns out, however, that equating  both sums of  Betti numbers 
 directly is quite intricate (see
 \cite{LiSc13} for the rank $2$ case). So, in Theorem
 \ref{thm:main2}, we prove a
 stronger statement instead: the $i^\mathrm{th}$ Betti number of
 the real part of the moduli space is equal to the $i^\mathrm{th}$ ascending diagonal sum of the Hodge
diamond of its complex part  (see Figure \ref{Hodge_diamond}).
This implies maximality since the moduli spaces in question have  torsion-free integral
cohomology \cite{AtBo83}.

The idea of deducing maximality from this stronger statement
comes from the old empirical observation that, although few maximal
real algebraic varieties are currently known, a seemingly large proportion of
them, that we suggest to call 
\textit{Hodge-expressive}, in fact satisfy this stronger property.  
\begin{defi}
A non-singular projective real algebraic variety $X$ is called Hodge-expressive 
if $H^*(X;\Z)$ is torsion free and 
  \[
  P_t(\R X)\ =\ H_{(t,1)}(X).
  \]
\end{defi}
  
Clearly, a Hodge-expressive variety is maximal. As mentioned above, such
varieties already appear in the literature and, in an informal way,
one may say that Hodge-expressive real algebraic varieties constitute
the basic maximal real algebraic varieties. In Section
\ref{sec:discussion}, we provide a brief panorama of known
Hodge-expressive varieties. Usually, Hodge-expressivity of a real
algebraic variety is a consequence of a prior proof of its
maximality. Here, we propose to go the other way round: maximality is
established as a consequence of Hodge-expressivity. From what we know,
this is also the strategy used in the aforementioned unpublished construction by
Itenberg and Viro. 

Note that Hodge-expressive varieties satisfy $\chi=\sigma$ (in the terminology of \cite{Bru21}), since
  \[
   \chi(\R X)\ =\ P_{-1}(\R
   X)\ =\ H_{(-1,1)}(X)\ =\ \sigma(X),
   \]
where the last equality is the Hodge index Theorem.

\subsection{Moduli spaces of vector bundles over real algebraic curves }
Let $\sC$ be a non-singular   connected projective complex algebraic
curve
of genus $g\geq 1$.
Given two integer numbers  $r\geq 1$ and $d\in\Z$, we denote by
$\M_\sC(r,d)$ the moduli space of semistable holomorphic vector bundles of 
rank $r$ and degree $d$ over $\sC$.
When $r$ and $d$ are coprime, the space 
$\M_\sC(r,d)$ is a non-singular complex projective variety of
dimension $r^2(g-1)+1$. For the rest of this discussion, we assume
that $g\geq 2$. Using the number-theoretic approach to the cohomology
of moduli spaces of vector bundles developed by Harder and Narasimhan
in \cite{HaNa75}, Desale and Ramanan obtained in \cite{DeRa75} a
recursive formula for the rational Poincaré polynomial of
$\M_\sC(r,d)$. Then, using gauge theory, Atiyah and Bott proved in
\cite{AtBo83} that the integral cohomology of $\M_\sC(r,d)$ is
torsion-free and gave an alternative proof of the recursive formula
by Desale and Ramanan. Later on, Zagier obtained a closed
formula for $P_t(\M_\sC(r,d))$ in \cite{Zagier}, and Earl and
Kirwan used a finite-dimensional analogue of the Atiyah--Bott approach to prove a recursive formula
for the Hodge polynomial $H_{(x,y)}(\M_\sC(r,d))$ in \cite{EaKi00}. A
consequence of these formulas (either recursive or closed) is that 
Betti and Hodge numbers of $\M_\sC(r,d)$ depend only on the topological data $g$ and $r$ (and not $d$).

When $\sC$ is a non-singular real projective curve, the
anti-holomorphic involution $\tau:\sC\longrightarrow\sC$ induces a
real structure $\E\longmapsto \overline{\tau^*\E}$ on $\M_\sC(r,d)$,
turning the moduli space $\M_\sC(r,d)$ into a real algebraic variety
as well. In fact, if $r$ and $d$ are coprime and $\sC$ has real
points, the real locus of $\M_\sC(r,d)$ consists exactly of
isomorphism classes of geometrically stable real vector bundles of
rank $r$ and degree $d$ (\emph{cf.}~\cite{Sc12}), where geometric stability of a real bundle $(\sE,\tau)$ means that the underlying complex bundle $\sE$ is stable. 
In \cite{LiSc13}, Liu and Schaffhauser developed a real analogue of
the Atiyah--Bott approach to obtain a recursive
formula, as well as a closed formula, computing 
$P_t(\R\M_\sC(r,d))$. As in the complex
situation, the Betti numbers of $\R\M_\sC(r,d)$ 
are seen \textit{a
  posteriori}  not to depend on $d$
but only on $g, r$ and
$b_0(\R\sC)$. Note that
$b_0(\R \M_\sC(r,d))=2^{b_0(\R\sC)-1}$ by \cite{Sc12}.

Given a line bundle $\Lambda$ of degree $d$ on $\sC$, one may also
consider the moduli space $\M_\sC(r,\Lambda)$ of semistable
holomorphic vector bundles of  
rank $r$ with fixed determinant $\Lambda$. By definition, this space is a fiber of the
determinant map $\M_\sC(r,d)\longrightarrow \mathrm{Pic}^d(\sC)$.
When $r$ and $d$ are coprime, the space
$\M_\sC(r,\Lambda)$ is a non-singular complex projective variety of
dimension $(r^2-1)(g-1)$, which is real if both $\sC$ and $\Lambda$
are 
real. Both  polynomials $H_{(x,y)}(\M_\sC(r,\Lambda))$ and 
 $P_t(\R\M_\sC(r,\Lambda))$ can be deduced
from those of $\M_\sC(r,d)$, see Section \ref{sec:reduction}.

\medskip
The next statement is the main result of this note.

\begin{thm}\label{thm:main2}
Let $\sC$ be a maximal non-singular real projective curve of genus
$g\geq 1$. Let $r\geq 1$ and $d\in\Z$ be coprime integers, and let
$\Lambda$ be a real line bundle of degree $d$ on $\sC$. Then both
moduli spaces $\M_\sC(r,d)$ and $\M_\sC(r,\Lambda)$ are
Hodge-expressive.
  \end{thm}

In particular, these moduli spaces are maximal when $\sC$ is maximal.
The next proposition shows that the converse is true for the moduli spaces
$\M_\sC(r,d)$.

\begin{prop}\label{prop:iff}
Let $\sC$ be a non-singular real projective curve of genus
$g\geq 1$ and such that $\R\sC\neq\emptyset$. Let $r\geq 1$ and $d\in\Z$ be coprime integers. If $\sC$ is not maximal, then neither is $\M_\sC(r,d)$.
  \end{prop}
  
Note that $\M_\sC(r,\Lambda)$ is maximal when either $g=1$ or $r=1$ (provided
$\R\sC\neq\emptyset$), since in this case it is
reduced to a point. 
Whether $\M_\sC(r,\Lambda)$ is non-maximal when $\sC$
is non-maximal and  $g,r\geq 2$ remains an 
open question.
As noted in \cite{LiSc13}, the explicit formulas for the mod $2$ Betti
numbers of $\R\M_\sC(r,\Lambda)$ quickly become too complex to be
evaluated at $t=1$ (and similar remarks apply for the Hodge numbers of
$\M_\sC(r,\Lambda)$). It was however checked in \cite{LiSc13}, using a
computer, that this is indeed the case up until $r=6$. 
As is visible from Section \ref{sec:reduction},  our proof of Proposition
\ref{prop:iff} does not allow us to conclude that $\M_\sC(r,\Lambda)$
is not maximal when $\sC$ is not maximal and $g,r\geq 2$.

\subsection{Outline of the paper}

Section  \ref{sec:proofs} is devoted to the proof of Theorem
\ref{thm:main2} and Proposition \ref{prop:iff}.
In Section \ref{sec:discussion}, we give a panorama of known
maximal and Hodge-expressive real algebraic varieties.

\subsection*{Acknowledgments}
  We are grateful to Ilia Itenberg and anonymous referees
  for pointing out to us
  examples of maximal real algebraic varieties that we had initially missed.

\section{Hodge-expressivity for moduli spaces of vector bundles}\label{sec:proofs}

In this section we prove Theorem \ref{thm:main2} and Proposition \ref{prop:iff}.
We first treat the particular cases $g=1$ and $r=1$ in Sections \ref{sec:g=1}
and \ref{sec:r=1} respectively.
As a consequence, we prove Proposition \ref{prop:iff}
in Section \ref{sec:reduction} and show that it suffices to prove Theorem
\ref{thm:main2} for $\M_\sC(r,d)$. This is achieved in Sections
\ref{sec:hodge} and \ref{proof_main_thm}.

\subsection{The case $g=1$}\label{sec:g=1}
If $\sC$ is of genus $1$ and $\mathrm{gcd}(r, d)=1$, the determinant map
$\M_\sC(r,d)\longrightarrow\mathrm{Pic}^d(\sC)$ is an isomorphism of
algebraic varieties \cite{Tu93}. In particular, one gets an isomorphism
$\M_\sC(r,d)\simeq \sC$
after  a point $x_0\in \sC$ has been chosen.
In the case when  $\sC$ is real and $x_0\in\R\sC$, this isomorphism is
real by \cite{BiSc16}, proving Theorem~\ref{thm:main2} 
and Proposition~\ref{prop:iff}
in the case
$g=1$.

\medskip
From now on, we suppose that $g\geq 2$.

\subsection{The case $r=1$}\label{sec:r=1}
The case of moduli spaces of line bundles over $\sC$
is well known, since these moduli spaces are the Picard varieties
$\mathrm{Pic}^d(\sC)$. In particular one has
\[
H_{(x,y)}\big(\mathrm{Pic}^d(\sC)\big)\ =\ (1+x)^g\,(1+y)^g.
\]
On the other hand by \cite{GrHa83}, one has 
\[
P_t\big(\R\mathrm{Pic}^d(\sC)\big)\ =\ 2^{b_0(\R\sC)-1} P_t((\mathbb{S}^1)^g)\ =\ 2^{b_0(\R\sC)-1} (1+t)^g
\]
when  $\R\sC\neq\emptyset$.
Since $\sC$ is maximal if and only if $b_0(\R\sC)=g+1$, 
this proves Theorem~\ref{thm:main2} and Proposition~\ref{prop:iff} in
the case $r=1$.

\subsection{Preliminaries}\label{sec:reduction}
Recall that the groups $H^*(\M_\sC(r,d);\Z)$ and
$H^*(\M_\sC(r,\Lambda);\Z)$ are torsion-free by \cite{AtBo83}. Hence
it is enough to prove the two identities
\[
P_t\big(\R\M_\sC(r,d)\big)=H_{(t,1)}\big(\M_\sC(r,d)\big)
\qquad \mbox{and}\qquad
P_t\big(\R\M_\sC(r,\Lambda)\big)=H_{(t,1)}\big(\M_\sC(r,\Lambda)\big).
\]
Moreover when $\mathrm{gcd}(r, d)=1$, it is proved in \cite{EaKi00}
and in \cite{Baird_MA}, respectively, 
that
\[
H_{(x,y)}\big(\M_\sC(r,d)\big)\ =\ H_{(x,y)}\big(\mathrm{Pic}^d(\sC)\big)
H_{(x,y)}\big(\M_\sC(r,\Lambda)\big)
\]
and
\[
P_t\big(\R
\M_\sC(r,d)\big)\ =\ P_t\big(\R\mathrm{Pic}^d(\sC)\big)
P_t\big(\R\M_\sC(r,\Lambda)\big).
\]
In view of Section \ref{sec:r=1}, this gives
\begin{equation}\label{useful_eq}
H_{(t,1)}\big(\M_\sC(r,d)\big)\ =\ 2^g\,(1+t)^g
H_{(t,1)}\big(\M_\sC(r,\Lambda)\big),
\end{equation}
and
\[
P_t\big(\R
\M_\sC(r,d)\big)\ =\ 2^{b_0(\R\sC)-1} (1+t)^gP_t\big(\R\M_\sC(r,\Lambda)\big).
\]
In particular, applying the Smith--Thom inequality to
$\M_\sC(r,\Lambda)$, we get:
\begin{align*}
\sum_{i\ge 0}b_i\big(\R
\M_\sC(r,d)\big)&\ =\ P_1\big(\R
\M_\sC(r,d)\big)
\\ &\ =\ 2^{b_0(\R\sC)-1+g}P_1\big(\R\M_\sC(r,\Lambda)\big)
\\ &\ \le\ 2^{2g} P_{1}\big(\M_\sC(r,\Lambda)\big)\ =\ \sum_{i\ge 0}b_i\big(\M_\sC(r,d)\big)
\end{align*}
and equality can only hold here if $b_0(\R\sC)=g+1$, that is to say if
$\sC$ is maximal. This proves Proposition~\ref{prop:iff}.\hfill\qed

\medskip
Next, when $\sC$ is maximal, we see from \eqref{useful_eq} that  $\M_\sC(r,d)$ is
Hodge-expressive if and only if $\M_\sC(r,\Lambda)$ is
Hodge-expressive. To prove Theorem
\ref{thm:main2}, it is therefore sufficient to  prove
\begin{equation}\label{specialization}
P_t\big(\R\M_\sC(r,d)\big) = H_{(t,1)}\big(\M_\sC(r,d)\big).
\end{equation}
 In order to do so, we will show in the next two sections,
using results from \cite{EaKi00,LiSc13}, that the polynomials
$H_{(t,1)}(\M_\sC(r,d))$ and $P_t(\M_\sC(r,d))$
satisfy the same recursion relation.

\subsection{Hodge numbers of  $\M_\sC(r,d)$}\label{sec:hodge}

Here we recast the computation from \cite{AtBo83,EaKi00} of  Poincaré and Hodge
polynomials of the 
moduli spaces $\M_\sC(r,d)$.
For all $(r,d)\in\Z_{>0}\times \Z$, we denote by $Bun_\sC(r,d)$ the
moduli stack of all vector bundles of rank $r$ and degree $d$ on
$\sC$. It contains, as an open substack, the moduli stack of
semistable vector bundles of rank $r$ and degree $d$ on $\sC$, which
we denote by $Bun^{\ss}_\sC(r,d)$. If we fix a $C^\infty$ complex vector bundle $E$ of rank $r$ and degree $d$ over $\sC$, there is an isomorphism of stacks $$Bun_\sC(r,d)\ \simeq\ [\cA_E/\cG_E]\, ,$$ where $\cA_E$ is the set of Dolbeault operators (\emph{i.e.} holomorphic structures) on $E$ and $\cG_E$ is the automorphism group of $E$ (\emph{i.e.} the gauge group). In particular, we can think of the cohomology of the stack $Bun_\sC(r,d)$ in the sense of \cite{Behrend} simply as the $\cG_E$-equivariant cohomology of the affine space $\cA_E$. By \cite{AtBo83}, the integral
cohomology of $Bun_\sC(r,d)$ is torsion-free. Moreover, its Poincaré
series does not depend on $d$ and is given by  
\[
P_t\big(Bun_\sC(r,d)\big)\ =\ \frac{(1+t)^{2g}}{1-t^2}\,\prod_{i=2}^r\frac{(1+t^{2i-1})^{2g}}{(1-t^{2i-2})(1-t^{2i})}\,.
\]
By \cite{HaNa75},  an algebraic vector bundle $\sE$ over $\sC$
admits a unique filtration, called the \textit{Harder--Narasimhan filtration},
\[
0=\sE_0 \subset \sE_1 \subset\, \cdots\, \subset \sE_\ell = \sE
\]
such that:
\begin{itemize}
\item for all $i\in\{1,\,\ldots\,,\ell\}$, the vector bundle $\sE_i/\sE_{i-1}$ is semistable;
\item setting $r_i := \mathrm{rk}(\sE_i/\sE_{i-1})$ and $d_i:=\mathrm{deg}(\sE_i/\sE_{i-1})$, we have $$\frac{d_1}{r_1} >\, \cdots\, > \frac{d_l}{r_l}\,.$$
\end{itemize}
The topological invariants $(r_i,d_i)_{1\leq i
  \leq\ell}$ of the successive quotients $\sE_i/\sE_{i-1}$ constitute the
\textit{Harder--Narasimhan type} (or HN type) of the vector bundle
$\sE$. If $\sE$ is  of rank $r$ and degree $d$, then its HN type is subject to the following constraints:

 \begin{enumerate}
  \item $d_1+\,\cdots\,+d_l=d$,
  \item $r_1+\,\cdots\,+r_l=r$,
  \item $\frac{d_1}{r_1} >\, \cdots\, > \frac{d_l}{r_l}\,.$
\end{enumerate}
By definition, the vector bundle $\sE$ is semistable if and only if
it is of HN type $(r,d)$, which we denote by $\mu_{\ss}$.
We denote by 
$I_{r,d}$ the set of all possible tuples of integers
$\mu=(r_i,d_i)_{1\leq i \leq\ell_\mu}$ satisfying properties (1), (2), (3) above.
Note that this set is infinite as soon as $r\geq 2$. For instance:
\[
I_{2,d}=\left\{(2,d), (1,k,1,d-k) \ |\ k > \frac{d}{2}
\right\}. 
\]

For all $\mu\in I_{r,d}$, there exists an algebraic substack
$Bun_\sC(\mu)$ of $Bun_\sC(r,d)$, parameterizing vector bundles of
type $\mu$, with
 $Bun_\sC(\mu_{\ss})=Bun_\sC^{\ss}(r,d)$.
The codimension of $Bun_\sC(\mu)$ is finite and, for all
$\mu\neq\mu_{\ss}$, it is equal to 
\[
d_\mu\ :=\ \sum_{1\leqslant i<j\leqslant l(\mu)}\big(r_j d_i - r_i d_j
+ r_i r_j (g - 1)\big).
\]
Moreover, 
there is an isomorphism of $\Q$-vector spaces
 \[
H^*\big(Bun_\sC(\mu);\mathbb{Q}\big)\ \simeq\ \bigotimes_{i=1}^{\ell_\mu}
H^*\big(Bun_\sC^{\ss}(r_i,d_i);\mathbb{Q}\big).
\]
In particular,
we have, for all $\mu\in I_{r,d}$,
\[
P_t\big(Bun_\sC(\mu)\big)\ =\ \prod_{i=1}^{\ell_\mu}
P_t\big(Bun_\sC^{\ss}(r_i,d_i)\big)\,.
\]
Finally, the stratification of the moduli stack $Bun_\sC(r,d)$ by the
substacks $Bun_\sC(\mu)$
is \textit{perfect} in the sense that the associated Gysin long exact
sequence breaks up into short exact sequences \cite{AtBo83}. This  implies that 
\[
P_t\big(Bun_\sC(r,d)\big)\ =\ \sum_{\mu\in I_{r,d}} t^{2d_\mu}\,
P_t\big(Bun_\sC(\mu)\big),
\]
yielding the recursive formula 
\[
P_t\big(Bun_\sC^{\ss}(r,d)\big)\ =\ P_t\big(Bun_\sC(r,d)\big)\ -
\sum_{\mu\neq\mu_{\ss}} t^{2d_\mu}\, \prod_{i=1}^{\ell_\mu}
P_t\big(Bun_\sC^{\ss}(r_i,d_i)\big)\,.
\]
When $r$ and $d$ are coprime, all semistable vector bundles of rank
$r$ and degree $d$ are stable.  Atiyah and Bott have shown that, in
this case, the Poincaré polynomial of  $\M_\sC(r,d)$ is
related to the Poincaré series of  $Bun_\sC^{\ss}(r,d)$
via the identity 
\[
P_t\big(\M_\sC(r,d)\big)\ =\ (1-t^2)\, P_t\big(Bun_\sC^{\ss}(r,d)\big)\,.
\]
Altogether, setting
$Q^{\C}_t(r,d)= (1-t^2)P_t\big(Bun_\sC^{\ss}(r,d)\big)$, one obtains the recursive formula
\begin{equation}\label{rec_for_Poincare_complex_mod_space}
Q^{\C}_t(r,d)\ =\ (1+t)^{2g}\prod_{i=2}^r\frac{(1+t^{2i-1})^{2g}} 
{(1-t^{2i-2})(1-t^{2i})}\ 
- \negmedspace\sum_{\mu\in I_{r,d}\setminus\{ \mu_{\ss}\}} \frac{t^{2d_\mu}}{(1-t^2)^{\ell_\mu-1}}\,
\prod_{i=1}^{\ell_\mu} Q^{\C}_t(r_i,d_i)\, ,
\end{equation} 
expressing all  polynomials $Q^{\C}_t(r,d)$
in terms of the initial term
\[
Q^{\C}_t(1,d)\ =\ P_t\big(\M_\sC(1,d)\big)\ =\ P_t\big(\mathrm{Pic}^d(\sC)\big)\ =\ (1+t)^{2g}\,.
\]

\medskip

In \cite{EaKi00}, Earl and Kirwan obtained a similar
recursive formula  for the Hodge polynomial of $\M_\sC(r,d)$ when $r$
and $d$ are coprime. Namely, setting $Q_{(x,y)}(r,d) =
H_{(x,y)}(\M_\sC(r,d))$
and using the construction of $\M_\sC(r,d)$ as a GIT quotient, they proved that 
\begin{equation}\label{eq:H}
  \begin{aligned}
      Q_{(x,y)}(r,d)\ =&\ (1-xy)\,G_{(x,y)}(r,d)\
-\negmedspace\sum_{\mu\in I_{r,d}\setminus\{ \mu_{\ss}\}} \frac{(xy)^{d_\mu}}{(1-xy)^{\ell_\mu-1}}
\, \prod_{i=1}^{\ell_\mu}
Q_{(x,y)}(r_i,d_i)\,,
  \end{aligned}
\end{equation}
expressing recursively all polynomials $Q_{(x,y)}(r,d)$
in terms of the power series 
\[
G_{(x,y)}(r,d):=\frac{(1+x)^g\,(1+y)^g}{1-xy}\,\prod_{i=2}^r\frac{(1+x^{i-1}y^i)^g\,
  (1+x^iy^{i-1})^g}{(1-x^{i-1}y^{i-1})(1-x^iy^i)}
\]
and the initial term
\[
Q_{(x,y)}(1,d)\ =\ (1-xy)\,G_{(x,y)}(1,d)\ =\ (1+x)^{g}(1+y)^g\ =\ H_{(x,y)}\big(\mathrm{Pic}^d(\sC)\big)\,.
\]

Plugging $x=y=t$ in equation~\eqref{eq:H}, one finds again the
recursion \eqref{rec_for_Poincare_complex_mod_space}. And plugging $x=t$ and $y=1$ in equation~\eqref{eq:H}, one obtains the following recursion:
\begin{equation}\label{eq:H2}
  \begin{aligned}
    Q_{(t,1)}(r,d)\ =&\
2^g(1+y)^g\,\prod_{i=2}^r\frac{(1+t^{i-1})^g\, (1+t^i)^g}{(1-t^{i-1})(1-t^i)} \
-\negmedspace\sum_{\mu\in I_{r,d}\setminus\{ \mu_{\ss}\}} \frac{t^{d_\mu}}{(1-t)^{\ell_\mu-1}}
\, \prod_{i=1}^{\ell_\mu}
Q_{(t,1)}(r_i,d_i)\,,
  \end{aligned}
\end{equation}
which will be useful to prove Theorem \ref{thm:main2}.

\begin{rem}
Although this is not how Earl and Kirwan prove their result, the power
series $G_{(x,y)}(r,d)$ 
can be interpreted as the Hodge series of the moduli stack $Bun_\sC(r,d)$ \cite{Te98,BeDi07}.
\end{rem}

\subsection{Betti numbers of $\R\M_\sC(r,d)$}\label{proof_main_thm}

When the curve $\sC$ is defined over the reals, so is the moduli stack
$Bun_\sC(r,d)$.
Here, we recast the computation of the mod $2$ Poincaré
polynomial of  $\R \M_\sC(r,d)$ from \cite{LiSc13} .
We denote by $n:=b_0(\R\sC)$ the number of connected components of
$\R\sC$,
and we assume that $n>0$. The stack $\R Bun_\sC(r,d)$ is a disjoint
union $$\R Bun_\sC(r,d)\ = \ \bigsqcup_s Bun_\sC^\R(r,d,s)$$ of
$2^{n-1}$ connected components, each of which is indexed by the
possible real invariants of a real vector bundle $(\sE,\tau)$ of rank
$r$ and degree $d$, namely the first Stiefel--Whitney
class $$s=(s_1,\,\ldots\,,s_{n})\in(\Z/2\Z)^{n}$$ of the vector bundle
$\R\sE$ over $\R\sC=\sqcup_{i=1}^{n} \mathbb{S}^1$, subject to the condition
$s_1+\,\cdots\,+s_n=d\mod{2}$. The substack $Bun_\sC^\R(r,d,s)$ is the
stack of all real vector bundles of rank $r$, degree $d$ and real type
$s$. Fixing  a $C^\infty$ real vector bundle $(E,\tau)$ with these
invariants $(r,d,s)$, we get an isomorphism of
stacks $$Bun_\sC^\R(r,d,s)\ \simeq\ [\cA_E^\tau/\cG_E^\tau]\, ,$$
where $\cA^\tau\subset\cA$ is the set of all $\tau$-fixed Dolbeault
operators (for the real structure on $\cA_E$ induced by the real
structure of $E$), and $\cG_E^\tau\subset\cG_E$ is the group of gauge
transformations of $E$ that commute with $\tau$. It turns out that the
mod $2$ Poincaré series of $Bun_\sC^\R(r,d,s)$ is independent of $d$ and
$s$ \cite{LiSc13}. As a consequence, the space
$\R Bun_\sC(r,d)$ has mod 2 Poincaré series
\[
  \ P_t\big(\R
  Bun_\sC(r,d)\big)\ =\ 2^{n-1}\,\frac{(1+t)^{g}}{1-t}\,\prod_{i=2}^r\frac{(1+t^{2i-1})^{g+1-n}\,(1+t^{i-1})^{n-1}\,(1+t^i)^{n-1}}{(1-t^{i-1})(1-t^{i})}\,.
  \]

Because of the uniqueness of the destabilizing sub-bundle, semistability
over $\R$ is equivalent to semistability over $\C$ and all sub-bundles
$(\sE_i)_{1\leq i\leq \ell}$ in the Harder--Narasimhan filtration of a
real algebraic vector bundle $\sE$ over $\sC$ are real. So the
substack $Bun_\sC(\mu)$ is defined over $\R$ and $\R Bun_\sC(\mu)$ is the stack of all real vector bundles of HN type $\mu$. Moreover, for all $\mu\in I_{r,d}$,
there is
an isomorphism of $(\Z/2\Z)$-vector spaces
\[
H^*\big(\R Bun_\sC(\mu);\ZZ\big)\ \simeq\ \bigotimes_{i=1}^{\ell_\mu}
H^*\big(\R Bun_\sC^{\ss}(r_i,d_i);\ZZ\big).
\]
In particular, we have, for all $\mu\in I_{r,d}$,
\[
P_t\big(\R Bun_\sC(\mu)\big)\ =\ \prod_{i=1}^{\ell_\mu}
P_t\big(\R Bun_\sC^{\ss}(r_i,d_i)\big)\,.
\]
Finally, the substacks $(\R Bun_\sC(\mu))_{\mu\in I_{r,d}}$ form a
stratification of $\R Bun_\sC(r,d)$ which is perfect over $\Z/2\Z$,
so 
\[
P_t\big(\R Bun_\sC(r,d)\big)\ =\ \sum_{\mu\in I_{r,d}} t^{d_\mu}\,
P_t\big(\R Bun_\sC(\mu)\big)
\]
yielding the recursive formula 
\[
 P_t\big(\R Bun_\sC^{\ss}(r,d)\big)\ =\ P_t\big(\R Bun_\sC(r,d)\big)\
  - \negmedspace\sum_{\mu\in I_{r,d}\setminus\{\mu_{\ss}\}} t^{d_\mu}\,
  \prod_{i=1}^{\ell_\mu}P_t\big(\R Bun_\sC^{\ss}(r_i,d_i)\big).
  \]
When $r$ and $d$ are coprime, one has the identity 
\[
P_t\big(\R \M_\sC(r,d)\big)\ =\ (1-t)\, P_t\big(\R
Bun_\sC^{\ss}(r,d)\big)\,.
\]
Altogether, setting
$Q^\R_t(n,r,d)= (1-t)P_t\big(\R Bun_\sC^{\ss}(r,d)\big)$, one obtains the recursive formula
\begin{equation}\label{rec_for_Poincare_real_mod_space}
  \begin{aligned}
  Q^\R_t(n,r,d)\ =\ \ &
  2^{n-1}\,(1+t)^{g}\,\prod_{i=2}^r\frac{(1+t^{2i-1})^{g+1-n}\,(1+t^{i-1})^{n-1}\,(1+t^i)^{n-1}}{(1-t^{i-1})(1-t^{i})}
 \\ & - \negmedspace\sum_{\mu\neq\mu_{\ss}} \frac{t^{d_\mu}}{(1-t)^{\ell_\mu-1}}\,
  \prod_{i=1}^{\ell_\mu} Q^\R_t(n,r_i,d_i), 
  \end{aligned}
\end{equation} 
expressing  all  polynomials $Q^\R_t(n,r,d)$
in terms of the initial term
\[
Q^\R_t(n,1,d)\
=\ P_{t}\big(\R\mathrm{Pic}^d(\sC)\big)\ =\ 2^{n-1}(1+t)^g\,.
\]

Thus it follows
from relations~\eqref{eq:H2}
and \eqref{rec_for_Poincare_real_mod_space}
that, for all $r$ and $d$, one has
\[
Q^\R_t(g+1,r,d)\ =\ Q_t(r,d).
\]
Since, when $d$ and $r$ are coprime,
\[
Q^\R_t(g+1,r,d)\ =\ P_t\big(\R\M_\sC(r,d)\big) \qquad\mbox{and}\qquad
Q_t(r,d)\ =\ H_{(t,1)}\big(\M_\sC(r,d)\big), 
\]
by results from \cite{EaKi00,LiSc13}, we see that relation~\eqref{specialization}, hence also Theorem~\ref{thm:main2}, is proved.\hfill\qed

\section{A brief panorama of Hodge-expressive varieties}\label{sec:discussion}

We list below a few basic examples of Hodge-expressive varieties and
discuss their relevance among maximal real algebraic varieties.

\subsection{Handy real algebraic varieties}\label{sec:handy}
Real projective spaces, \textit{i.e.} $\CP^n$ equipped with the
standard complex conjugation, constitute the simplest  
projective  real algebraic varieties. One checks easily that they are
all Hodge-expressive.

More generally, the standard complex conjugation in $\C$ induces a
real structure on the Grassmannian variety $\mathbf{Gr}(d,n)$
of $d$-planes in
$\C^n$. Schubert cells provide a stratification
of $\mathbf{Gr}(d,n)$ by real affine spaces \cite{MilSta74}, implying
that
for all $i\ge 0$
\[
b_{2i+1}\big(\mathbf{Gr}(d,n)\big)=0
\qquad\mbox{and}\qquad
b_i\big(\R \mathbf{Gr}(d,n)\big)=b_{2i}\big(\mathbf{Gr}(d,n)\big)
=h^{i,i}\big(\mathbf{Gr}(d,n)\big).
\]
Hence all Grassmannian varieties are Hodge-expressive.

Toric varieties constitute another generalization of projective
spaces. A toric variety carries a standard real structure inherited
once again from the standard complex conjugation in $\C^*$, and 
all these non-singular real projective toric varieties are Hodge-expressive
by \cite{BFMV06,Fr22}.

Grassmannian and toric varieties are particular cases of
{\it balanced varieties}, \textit{i.e.} non-singular projective complex
manifold $X$ such that $h^{i,j}(X)=0$ whenever $i\neq j$.
By definition, such a variety satisfies
$P_t(\C X)=H_{(t,t)}(X) =H_{(t^2,1)}(X)$. In particular if $X$ is
both balanced and Hodge-expressive, then 
$P_t(\C X)\ =\ P_{t^2}(\R X)$.

\medskip
Hodge-expressivity is preserved by certain elementary operations on real algebraic varieties, for instance taking products or projectivizing a vector bundle over the variety. This follows from the fact that both Hodge and
Betti numbers satisfy the same relations under these elementary operations:
on the one hand, the Hodge
polynomial extends to a motivic invariant from the Grothendieck ring
$K_0(Var_\C)$ 
of complex algebraic varieties to $\Z[x,y]$ \cite[Remark 5.56]{PetSte08}, and on the
other hand, the
 Künneth formula (or more generally the Leray--Hirsch Theorem,
 \cite[Theorem 4D.1]{Hat02}) implies that the Poincaré polynomial behaves analogously
in basic situations. More precisely, we observe the following:

\begin{itemize}
\item The product of Hodge-expressive real algebraic varieties $X_1, X_2,\,\ldots\,, X_k$,
equipped with the product real structure, is Hodge-expressive. Indeed,
one deduces from the multiplicativity of the Hodge polynomial and from
the Künneth formula that
\[
H_{(x,y)}\left(\prod_{i=1}^k X_i\right)=\prod_{i=1}^k H_{(x,y)}(X_i) ,
\quad\mbox{and}
\quad
P_{t}\left(\prod_{i=1}^k \R X_i\right)=\prod_{i=1}^k
P_{t}(\R X_i).
\]
\item Similarly, if $E$ is a real algebraic vector bundle of rank $r$ over a Hodge-expressive real
algebraic variety $X$,  the real projective bundle $\mathbb P(E)$
is also Hodge-expressive. By the scissor relations on  $K_0(Var_\C)$, we get that
\[
H_{(x,y)}(\mathbb P(E))=H_{(x,y)}(X)H_{(x,y)}(\CP^{r-1}).
\]
And since the first Chern class (resp. the first
Stiefel--Whitney class) of the tautological line
bundle on $\mathbb P(E)$ (resp. $\R\mathbb P(E)$) restricts to the
generator of $H^*(\CP^{r-1})$ (resp. $H^*(\RP^{r-1})$), we deduce from the Leray--Hirsch theorem that
$H^*(\mathbb P(E);\Z)$ is torsion-free if $X$ is, and that
\[
P_{t}(\R\mathbb P(E))=P_{t}(X)P_{t}(\RP^{r-1}).
\]
\item
  We deduce in particular from the last item that the blow-up $\widetilde X$ of a Hodge-expressive
variety $X$ along a
Hodge-expressive subvariety $Y$ is Hodge-expressive as well.
Using again the scissor relations on  $K_0(Var_\C)$, we have
\[
H_{(x,y)}(\widetilde X)=H_{(x,y)}(X) + (x^{r-1}y^{r-1} +
x^{r-2}y^{r-2}+\cdots +xy) H_{(x,y)}(Y),
\]
where $r$ is the codimension of $Y$ in $X$.
The computation of Betti numbers of $\widetilde X$ from
\cite[Chapter~4, Section~6]{GrHa78} carries over word for word to the real
part. In particular,  we obtain that $H^*(\widetilde X;\Z)$ is
torsion-free if $X$ and $Y$ are, and that
\[
P_{t}(\R\widetilde X)=P_{t}(\R X) + (t^{r-1} +
t^{r-2}+\cdots +t) P_{t}(\R Y).
\]
\end{itemize}

\smallskip

One may also consider symmetric powers of real algebraic varieties, or
more generally  quotients $X^k/\Gamma$ of a $k^\mathrm{th}$ product
$X^k$ by the canonical action of a
subgroup $\Gamma$ of the symmetric group $\mathfrak S_k$. The product
and the symmetric power 
correspond to the two extremal cases $\Gamma=\{\mathrm{Id}\}$ and
$\Gamma=\mathfrak S_k$, respectively. Any real structure on
$X$ extends canonically to  $X^k/\Gamma$, and
Franz proved in \cite{Fr18} that  $X^k/\Gamma$ is 
a maximal real algebraic variety as soon as $X$ is (note that the
converse does not hold in general, as shown by the second symmetric power
of $\CP^1$ equipped with the antipodal involution).
It would be interesting to investigate whether
$X^k/\Gamma$ is Hodge-expressive as soon as it is
non-singular and $X$ is Hodge-expressive.
When  $\Gamma=\mathfrak S_k$, the $k^\mathrm{th}$ symmetric power
$X^{[k]}$ of $X$
is projective and non-singular if and only if $X$ is a non-singular
projective curve.
In this case partial results regarding Hodge-expressivity can
be obtained, as we shall see in Section \ref{sec:curves} below.

\subsection{Curves}\label{sec:curves}
 Let $\sC$ be a non-singular real projective curve of genus $g$. In
 this case the Smith--Thom inequality  reduces to
 the Harnack--Klein inequality
 \[
 b_0(\R \sC)\ \leq\ g+1,
 \]
 with equality if and only if $\sC$ is maximal.
Since
\[
H_{(x,y)}(\sC) = 1+g(x+y) + xy,
\]
and
\[
P_t(\R\sC)\ =\  b_0(\R \sC) (1+t)\le (g+1)(1+t)\ =\ H_{(t,1)}(\sC),
\]
the curve $\sC$ is Hodge-expressive if and only if it is maximal. 

\medskip
By the Abel--Jacobi Theorem, we have
$\mathrm{Pic}^d(\sC) \simeq \mathrm{Jac}(\sC)$ for all $d\in\Z$, so
\[
H_{(x,y)}\big(\mathrm{Pic}^d(\sC)\big)\ =\ (1+x)^g\,(1+y)^g.
\]
By \cite{GrHa83}, if $b_0(\R \sC)> 0$ one has
\[
P_t\big(\R\mathrm{Pic}^d(\sC)\big)\ =\ 2^{b_0(\R \sC)-1}\,(1+t)^g\ \leq\ 2^{g}\,(1+t)^g\ =\ H_{(t,1)}\big(\mathrm{Pic}^d(\sC)\big).
\]
Since $\mathrm{Pic}^d(\sC)$ has torsion-free integral cohomology, we
deduce that
\begin{equation}\label{equiv_PicC_and_C}
\mathrm{Pic}^d(\sC) \mbox{ Hodge-expressive}
\Longleftrightarrow \mathrm{Pic}^d(\sC) \mbox{ maximal }
\Longleftrightarrow \sC \mbox{ maximal}
\end{equation}
as soon as  $b_0(\R \sC)> 0$.

\begin{rem}\label{Picard_var_when_RC_is_empty}
When $\R\sC=\emptyset$, one has $\R\mathrm{Pic}^d(\sC)=\emptyset$ if
$g$ and $d$ are both odd \cite{GrHa83}. And if $g$ is odd and $d$ is
even, the Poincaré polynomial of $\R\mathrm{Pic}^d(\sC)$ is
$2(1+t)^g$, so its value at $t=1$ is strictly smaller than $2^{2g}$,
unless $g=1$ (in which case $\mathrm{Pic}^{2d}(\sC)$ is maximal even
though $\R\sC=\emptyset$). Finally, if $g$ is even, the Poincaré
polynomial of $\R\mathrm{Pic}^d(\sC)$ is $(1+t)^g$, so its value at
$t=1$ is strictly smaller than $2^{2g}$, unless $g=0$. 
\end{rem}

From the utmost right equivalence in \eqref{equiv_PicC_and_C} and the discussion in Section \ref{sec:handy}, one deduces that the $k^\mathrm{th}$ symmetric power $\sC^{[k]}$
of $\sC$ is Hodge-expressive if $\sC$ is maximal and
$k\ge 2g-1$, which refines \cite[Theorem 3.1]{BisDme17}.
Indeed, the Riemann--Roch and Abel--Jacobi theorems
imply that $\sC^{[k]}$ can be
expressed as the projectivization of a real algebraic vector bundle
over $Jac(\sC)$, see \cite[Chapter IV]{ACGH85}. 
Furthermore, from the explicit computation of Betti numbers of $\R\sC^{[k]}$ for
$k=2,3$ in \cite{BisDme17}, one also sees that  $\sC^{[k]}$ is
Hodge-expressive when $\sC$ is maximal and $k=2,3$ (see for example
\cite{LS11} for Hodge numbers of $\sC^{[k]}$).
It would be interesting to investigate whether $\sC^{[k]}$ is
Hodge-expressive as soon as $\sC$ is.

\medskip

Hodge-expressivity may also be used to determine
``elementary embeddings'' of real algebraic curves in 
 real algebraic surfaces.
 As an example, let us consider
 Harnack curves in the real projective plane. They
were originally
 constructed by Harnack in \cite{Har}, and constituted the first
 family of maximal curves of arbitrary degree in $\RP^2$.
 Since then, Harnack curves have also been obtained
 as elementary applications of
 a large proportion of construction methods of maximal real algebraic curves
(\textit{e.g.} \cite{IV2}). Surprisingly, they also appeared in several
contexts other than pure real algebraic geometry (\textit{e.g.} \cite{KenOko06}). 
To summarize informally, a maximal real algebraic curve in $\RP^2$  is a Harnack curve, except if it  has a
good reason not to be so.

The topological type of the pair $(\RP^2, \R\sC)$ for a Harnack curve
$\sC$ in $\RP^2$ is depicted in Figure~\ref{fig:harnack} below.

\begin{figure}[h!]
\centering
\begin{tabular}{ccc} 
  \includegraphics[height=1cm, angle=0]{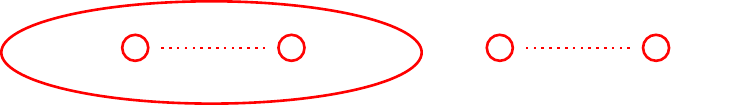}
  \put(-172,-5){$\smash{\underbrace{\qquad\qquad\quad }_{\frac{(k-1)(k-2)}{2}}}$}
  \put(-74,-5){$\smash{\underbrace{\qquad\qquad\quad }_{\frac{3k(k-1)}{2}}}$}
  & \hspace{5ex} &
\includegraphics[height=1cm, angle=0]{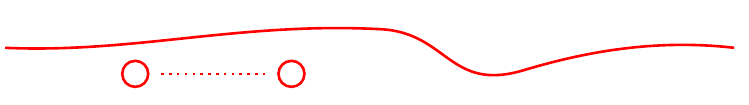}
  \put(-172,-5){$\smash{\underbrace{\qquad\qquad\quad }_{\frac{(d-1)(d-2)}{2}}}$}
\\ \\ \\ (a)\, $d=2k$ && (b)\, $d=2k+1$
\end{tabular}
\caption{Harnack curves of degree $d$ in $\RP^2$, up to isotopy and in
an affine chart.}
\label{fig:harnack}
\end{figure} 

\noindent Suppose that a Harnack curve $\sC$ of even degree $2k$ in $\RP^2$
is given by the real equation $P(x,y,z)=0$, where the real polynomial
$P(x,y,z)$ is positive on the non-orientable
connected component of $\RP^2\setminus\R\sC$. Then, on the one hand, one easily sees from
Figure \ref{fig:harnack} that the real algebraic surface $X$ with
equation
$w^2-P(x,y,z)$ in the weighted projective space $\CP(1,1,1,k)$
satisfies
\[
b_0(\R X)\ =\ b_2(\R X)\ =\ \frac{(k-1)(k-2)}{2}+1 \qquad\mbox{and}\qquad
b_1(\R X)\ =\ 3k(k-1)+2.
\]
On the other hand, the Hodge diamond  of $X$ is (see for example
\cite[Section 2.5]{DK}) 
 \[
 \begin{array}{ccccc}
 &  &1&&
\\ &  0 & & 0&
\\ \frac{(k-1)(k-2)}{2} &  &3k(k-1)+2 & &\frac{(k-1)(k-2)}{2}
\\  & 0 & & 0&
 \\ & &1&&
 \end{array}.
 \]
Note that $H^*(X;\Z)$ is torsion-free by the Lefschetz hyperplane section
Theorem. So the basic nature of 
a Harnack curve $\sC$ 
may be reflected in the fact
that $X$ is Hodge-expressive.

 \subsection{Surfaces}
 We refer to \cite[Section 3]{DK} and \cite[Chapter 4]{Mang17} for
 more details and references
 regarding the various classifications of real algebraic
 surfaces discussed in this section.
 
 It follows from the
 classifications by Comessati that all maximal real Abelian surfaces are Hodge-expressive. For
  examples, the Hodge diamond of an Abelian surface is
   \[
 \begin{array}{ccccc}
 &  &1&&
\\ &  2 & & 2&
\\ 1 &  &4 & &1
\\  & 2 & & 2&
 \\ & &1&&
 \end{array},
 \]
and such a real surface $X$ is maximal if and only if $\R X$ is the
disjoint union of four tori $\mathbb{S}^1\times \mathbb{S}^1$.

\medskip
There exist maximal real algebraic
$K3$-surfaces that are not Hodge-expressive. As shown by Kharlamov \cite{Kha76}, the
three possible topological types for the real part of a
maximal real $K3$-surface are: the disjoint union of a sphere $S^2$ and a
surface of genus $10$,  the  disjoint union of five spheres with a
surface of genus $6$, and the disjoint  union of nine spheres with a
surface of genus $2$.
Since a $K3$ surface has torsion-free integral cohomology and has
the following Hodge diamond
   \[
 \begin{array}{ccccc}
 &  &1&&
\\ &  0 & & 0&
\\ 1 &  &20 & &1
\\  & 0 & & 0&
 \\ & &1&&
 \end{array},
 \]
only the first topological type corresponds to a Hodge-expressive
variety.

\medskip

Note that  there exist  maximal real algebraic Enriques surfaces
satisfying
  \[
  b_i(\R X)\ >\ \sum_{j\geq 0} h^{i,j}(X)\quad \forall i\in\{0,1,2\},
  \]
see \cite{DK,Mang17}.  This shows the necessity of
  the torsion-freeness assumption in the implication 
    \[
 \Big(P_t(\R X)=H_{(t,1)}(X)\Big)\Longrightarrow X\mbox{ maximal}.
    \]
    
\medskip

The first family of maximal real algebraic surfaces of arbitrary degree in
$\CP^3$ was constructed by Viro in  \cite{Vir79}, 
 generalizing Harnack's construction.
A non-singular hypersurface in $\CP^n$ has torsion-free
integral cohomology by the Lefschetz hyperplane section Theorem,
and it turns out that 
all
surfaces in Viro's family are 
Hodge-expressive.

 \subsection{Higher dimension}
As mentioned earlier, our knowledge of maximal real
algebraic varieties of dimension at least 3 is quite restricted.
By \cite{Kra09}, there exists a unique topological type of maximal real
cubic 3-folds in $\CP^4$, and one checks that it is Hodge-expressive.
By \cite{FinKha10}, there exist three topological types of maximal real
cubic 4-folds in $\CP^5$,  only one of them being Hodge-expressive.

All maximal real algebraic projective hypersurfaces from
 Itenberg and Viro's unpublished construction are Hodge-expressive.
 As mentioned earlier,
 this is even a crucial point in their argumentation,
  since they
  prove that the real algebraic hypersurfaces that they construct are maximal
  by showing that they are Hodge-expressive.
 Note that the Itenberg--Viro's construction uses primitive combinatorial
 patchworking. 
Renaudineau and Shaw recently proved \cite{RenSha18}
that all maximal real algebraic hypersurfaces of a non-singular compact toric
variety that are obtained by primitive combinatorial patchworking are
Hodge-expressive, thus confirming a long standing conjecture of Itenberg's.

\smallskip

We are not aware of any general construction of maximal real 
hypersurfaces in $\CP^n$ that are not Hodge-expressive.


\end{document}